\documentclass[conference]{IEEEtran}
\IEEEoverridecommandlockouts
\usepackage{cite}
\usepackage{amsmath,amssymb,amsfonts,mathtools}
\usepackage{algorithm}
\usepackage{algorithmic}
\usepackage{graphicx}
\usepackage{textcomp}
\usepackage{subcaption}
\usepackage{hyperref}
\usepackage{xcolor}
\usepackage{enumitem}
\usepackage{caption}
\usepackage{epstopdf}
\mathtoolsset{showonlyrefs=true}

\newtheorem{definition}{Definition}
\newtheorem{assumption}{Assumption}
\newtheorem{theorem}{Theorem}

\graphicspath{{Results/}{Results/MRI/}{Results/Phantom/}}

\def\BibTeX{{\rm B\kern-.05em{\sc i\kern-.025em b}\kern-.08em
    T\kern-.1667em\lower.7ex\hbox{E}\kern-.125emX}}
\begin{document}

\title{A scaled, inexact and adaptive Fast Iterative Soft-Thresholding Algorithm
for convex image restoration
\thanks{The authors acknowledge the support provided by the French IEA CNRS VaMOS grant and by the Italian INdAM GNCS research group.}
}

\author{\IEEEauthorblockN{Luca Calatroni}
\IEEEauthorblockA{\textit{CNRS, UCA, INRIA} \\
\textit{2000 Route des Lucioles, 06903}\\
Sophia-Antipolis, France \\
\href{mailto:calatroni@i3s.unice.fr}{calatroni@i3s.unice.fr}}
\and
\IEEEauthorblockN{Simone Rebegoldi}
\IEEEauthorblockA{\textit{Dipartimento di Ingegneria Industriale} \\
\textit{Università degli studi di Firenze}\\
Firenze, Italy \\
\href{mailto:simone.rebegoldi@unifi.it}{simone.rebegoldi@unifi.it}}
}

\newcommand{\R}{\mathbb{R}}
\newcommand{\N}{\mathbb{N}}

\newtheorem{Theorem}{Theorem}[section]
\newtheorem{Def}{Definition}[section]
\newtheorem{Cor}{Corollary}[section]
\newtheorem{Lemma}{Lemma}[section]
\newtheorem{Pro}{Proposition}[section]
\newtheorem{assumptions}{Assumption}
\newtheorem{Remark}{Remark}[section]

\newenvironment{ifelse}{%
  \begin{list}{}{%
      \setlength{\topsep}{0pt}\setlength{\parskip}{0pt}
      \setlength{\partopsep}{0pt}\setlength{\itemsep}{0pt}}}%
  {\end{list}}

\makeatletter
  \newcounter{listi}
  \newcounter{listii}[listi]
  \newcounter{listiii}[listii]
  \newcounter{listiv}[listiii]
  \newcounter{listv}[listiv]
  \newcounter{listvi}[listv]

  \gdef\listctr{list\romannumeral\the\@listdepth}\expandafter
  \newcommand{\steplabel}[1]{\hfil\textsc{Step\ #1.}}
\makeatother

\newenvironment{AlgorithmSteps}[1][1]{%
  \begin{list}{\csname label\listctr\endcsname}{%
      \usecounter{\listctr}
      \let\makelabel\steplabel
      \settowidth{\labelwidth}{\textsc{Step\ #1.}}%
      \setlength{\leftmargin}{\labelwidth}\addtolength{\leftmargin}{\labelsep}}}%
  {\end{list}}

\newcommand{\stepname}[1]{\emph{#1}.}
\def\proof{{\it Proof. }}
\def\x{ {\ve x}}

\def\y{ {\ve y}} \def\tildey{\tilde{\y}}
\def\u{ {\ve u}}\def\lamk{{\lambda^{(k)}}}
\def\w{{\ve w}}
\def\z{ {\ve z}}
\def\g{ {\ve g}}
\def\d{{\ve d}}\def\p{{\ve p}}
\def\xk{ \x^{(k)}}
\def\txk{ \tilde{x}^{(k)}}
\def\yk{ \y^{(k)}}
\def\tildeyk{ \tilde{\y}^{(k)}}
\def\uk{ \u^{(k)}}
\def\zk{ \z^{(k)}}\def\zkk{ \z^{(k+1)}}
\def\xs{\x^*}\def\barx{\bar{\x}}
\def\yk{ \y^{(k)}}
\def\gk{ \w^{(k)}}\def\dk{ \d^{(k)}}
\def\xkk{ \x^{(k+1)}}
\def\xkm{ \x^{(k-1)}}
\def\ukm{ \u^{(k-1)}}
\def\vkm{ \v^{(k-1)}}
\def\ek{{\epsilon_k}}
\def\lmin{\lambda_{min}}\def\lmax{\lambda_{max}}
\def\ykk{ \y^{(k+1)}}
\def\ukk{ \u^{(k+1)}}
\def\xtilde{\tilde\x}\def\kinN{{k\in\mathbb N}}
\def\diam{\mbox{diam}}
\def\dom{\mbox{dom}}
\def\R{\mathbb R}
\def\H{\mathcal H}
\def\frec{f^{rec}}
\def\flev{f^{lev}}
\def\erec{\epsilon^{rec}}
\def\endproof{\hfill$\square$\vspace{0.3cm}\\}
\def\kinN{{k\in\N}}\def\prox{{\mbox{prox}}}
\newcommand\D{{\mathcal D}}
\newcommand{\stl}{\alpha}
\newcommand{\proj}[1]{\p_{\stl,\gamma}({#1})} \newcommand{\projj}[1]{\p_{\stl,D}({#1})}
\newcommand{\KL}{\operatorname{KL}}
\newcommand{\HS}{\operatorname{HS}}
\newcommand{\projk}[2][k]{\p_{\stl_{#1},\gamma_{#1}}({#2})}\newcommand{\projjk}[2][k]{\p_{\stl_{#1},D_{#1}}({#2})}
\newcommand{\Q}{{Q_{\stl,\gamma}}}
\newcommand{\Qk}[1][k]{{Q_{\stl_{#1},\gamma_{#1}}}}\newcommand{\xkb}{\x^{(k-1)}}
\newcommand{\tk}{t_k}\newcommand{\tkk}{t_{k+1}}\newcommand{\tkm}{t_{k-1}} \newcommand{\ak}{\stl_k}\newcommand{\akk}{\stl_{k+1}}
\newcommand{\vk}{v_k}\newcommand{\vkk}{v_{k+1}}\newcommand{\Dk}{{D_k}}\newcommand{\Dkk}{{D_{k+1}}}
\newcommand{\gammak}{\gamma_k} \newcommand{\gammakk}{{\gamma_{k+1}}}
\newcommand{\thk}{\theta_k}\newcommand{\thkk}{\theta_{k+1}}\newcommand{\thkb}{\theta_{k-1}}
\newcommand{\sk}{s_k}\newcommand{\skk}{s_{k+1}}
\def\x{ {x}}\def\y{ {y}}
\def\xkk{ {\x^{(k+1)}}} \def\ykk{ {\y^{(k+1)}}}
\def\xk{ \x^{(k)} }\def\yk{ \y^{(k)} }
\def\xkm{ \x^{(k-1)} }
\def\prox{\mathrm{prox}} \def\arg{\mathrm{arg}}
\newcommand{\argmin}{\operatornamewithlimits{argmin}}
\def\dom{\mbox{dom}}
\newcommand{\ri}{\mbox{ri}}
\newcommand{\txtr}[1]{\textcolor{red}{#1}}
\newcommand{\txtb}[1]{\textcolor{blue}{#1}} 

\maketitle

\begin{abstract}
In this note, we consider a special instance of the scaled, inexact and adaptive generalised Fast Iterative Soft-Thresholding Algorithm (SAGE-FISTA) recently proposed in \cite{SAGE-FISTA} for the efficient solution of strongly convex composite optimisation problems. In particular, we address here the sole (non-strongly) convex optimisation scenario, which is frequently encountered in many imaging applications. The proposed inexact S-FISTA algorithm shows analogies to the variable metric and inexact version of FISTA studied in \cite{Bonettini2018a}, the main difference being the use of an adaptive (non-monotone) backtracking strategy allowing for the automatic adjustment of the algorithmic step-size along the iterations (see \cite{Scheinberg-2014,Calatroni-Chambolle-2019}). 
A quadratic convergence result in function values depending on the backtracking parameters and the upper and lower bounds on the spectrum of the variable metric operators is given. Experimental results on TV image deblurring problems with Poisson noise are then reported for numerical validation, showing improved computational efficiency and precision.
\end{abstract}

\medskip


\section{Introduction}

The use of forward-backward (FB) algorithms is nowadays extremely popular in the context of variational imaging due to their easy applicability in many problems and their provable fast convergence when coupled with suitable  inertial updates, as in the case of the celebrated Fast Iterative Soft-Thresholding Algorithm (FISTA)  \cite{Beck-Teboulle-2009b}. Their successful application relies, in particular, on some practical assumptions on the composite energy minimization problem one aims to solve. First, a closed-form computation of the backward (proximal) step is desirable to avoid the use of inner solvers. Secondly, an accurate estimation of the `steepness' of the smooth component of the functional to minimize is required to provide a sufficiently meaningful forward update, thus avoiding an unnecessary large number of iterations till convergence. Whenever either (or both) of these two features is missing, the practical effectiveness of FB-type algorithms may be limited. To circumvent the former issue, recent approaches deal directly with the inexact calculation of the proximal point and provide appropriate conditions on the accuracy of these approximations which guarantee the same convergence properties of FISTA, see, e.g., \cite{Bonettini2018a,Bonettini2019}. In order to deal with the latter bottleneck, adaptive backtracking procedures favouring the adjustment of the algorithmic step-size along the iterations can be used, see, e.g.,  \cite{Scheinberg-2014,Calatroni-Chambolle-2019}.  Moreover, as previous studies in the context of smooth convex optimisation problems showed \cite{Lanteri-etal-2001}, FB algorithms may also benefit from suitable scaling approaches defined in terms of second-order (Newton-type) information. Such procedures have been shown to render particularly effective in the context of signal-dependent image reconstruction problems in a variety of applications ranging from biological to astronomical imaging, see e.g., \cite{Bertero2018}. Recently,  in \cite{SAGE-FISTA} the authors proposed a general inexact, scaled and adaptive FISTA-type algorithm  (named there SAGE-FISTA) designed for solving possibly strongly-convex composite problems and encompassing all the aforementioned approaches. Upon suitable conditions on the scaling updates and on the inexactness parameters sequence, accelerated convergence rates for the function values are there rigorously proved. 

\medskip

In this short note, we consider a particular case of SAGE-FISTA, that is we specify its definition, features and convergence guarantees when strong convexity is not explicitly taken into account, i.e. when the strong convexity parameter denoted by $\mu\geq 0$ in \cite{SAGE-FISTA} is set $\mu=0$. The resulting algorithm thus takes the form of a general convex inertial iterative scheme which shows analogies with adaptive \cite{Scheinberg-2014}, and inexact and variable metric \cite{Bonettini2018a} FISTA-type schemes previously considered in the literature. Upon suitable assumptions on the sequence of parameters describing inexactness and on the variable metric operators, we specify in Theorem \ref{th:convergence} the quadratic convergence result for the function values of the inexact S-FISTA iterates and report in Section \ref{sec:choice} a user's guide on how to verify these conditions in practice. The algorithm is finally validated on a classical image deblurring problem where Total Variation (TV) regularisation is combined with a Kullback-Leibler data fidelity and a non-negativity constraint, which is frequently encountered in the framework of image restoration problems under the presence of signal-dependent Poisson noise.

\section{Problem setting}
We start setting some notation. For a given Hilbert space ($\H$, $\langle \cdot, \cdot\rangle$), we denote by $\|\cdot\|$ the norm induced on $\H$ by $\langle \cdot, \cdot\rangle$. For any function $f:\mathcal{H}\rightarrow \R\cup\{+\infty\}$, we denote by $\dom(f):=\{\x\in\H: \ f(\x)< + \infty\}$ its domain.  For a continuously differentiable function $f:\mathcal{H}\rightarrow \mathbb{R}$, we then denote by $\mathbb{D}_f(x,y):=f(x)-f(y)-\langle \nabla f(y),x-y\rangle$ the Bregman distance of $f$ between the ponts $x, y\in\mathcal{H}$.
Denoting by $\mathcal{S}(\H)$ the set of linear, bounded and self-adjoint operators from $\H$ to $\H$, and by $\mathcal{I}\in \mathcal{S}(\mathcal{H})$ the identity operator on $\mathcal{H}$, we further recall the standard Loewner partial ordering relation defined  on $\mathcal{S}(\H)$, which,  for all $ D_1, D_2 \in \mathcal{S}(\H)$, reads:
$$
 D_1\preceq D_2 \Leftrightarrow \langle D_1 x, x \rangle \leq \langle D_2 x, x \rangle  \ \forall x \in \H \, .
$$
For $\eta_{inf},\eta_{sup}\in\mathbb{R}_{>0}
$ with $0<\eta_{inf}\leq \eta_{sup}$, we further introduce the following sets
\begin{align*}
\mathcal{D}_{\eta_{inf}}&:=\{D\in \mathcal{S}(\H): \ \eta_{inf} \mathcal{I} \preceq D \},\\
\mathcal{D}_{\eta_{inf}}^{\eta_{sup}}&:=\{D\in \mathcal{S}(\H): \ \eta_{inf} \mathcal{I} \preceq D \preceq \eta_{sup} \mathcal{I} \},
\end{align*}
and notice that, by definition, $\mathcal{D}_{\eta_{inf}}^{\eta_{sup}}\subseteq \mathcal{D}_{\eta_{inf}}$. If $D\in \mathcal{D}_{\eta_{inf}}$, then we have that
\begin{equation}\label{eq:inner_productD}
(x,y):=\langle Dx, y \rangle,
\end{equation}
defines an inner product on $\mathcal{H}$, and the notation $\|x\|_D:=\sqrt{\langle Dx, x\rangle}$ can be used to denote the $D-$norm induced by \eqref{eq:inner_productD}. By definition, we thus have that if $D\in \mathcal{D}_{\eta_{inf}}^{\eta_{sup}}$, the following inequality holds
\begin{equation}\label{ine_norm}
\eta_{inf} \| x\|^2 \leq  \|x\|_D^2 \leq \eta_{sup} \| x \|^2, \quad \forall \ x\in\H.
\end{equation}
Note that by \cite[Theorem 4.6.11]{Debnath-etal-1990}, we have that if $D\in\mathcal{D}_{\eta_{inf}}$, then $D$ is invertible and, if also $D\in\mathcal{D}_{\eta_{inf}}^{\eta_{sup}}$, then we have $D^{-1}\in\mathcal{D}_{1/\eta_{sup}}^{1/\eta_{inf}}$. 
For $\tau>0$, $D\in \mathcal{D}_{\eta_{inf}}^{\eta_{sup}}$ and a proper, convex, and lower semicontinous $g:\H\rightarrow \R\cup\{+\infty\}$, we further define the proximal operator of $g$ w.r.t. the metric induced by $D$ as 
$$
\operatorname{prox}_{g}^D(x)=\underset{z\in\mathcal{H}}{\operatorname{argmin}} \ g(z)+\frac{1}{2}\|z-x\|_D^2, \quad \forall \ x\in\mathcal{H}.
$$
Finally, we recall that for a given $Y\subseteq \mathcal{H}$ nonempty, closed convex set and $D\in \mathcal{D}_{\eta_{inf}}$, the projection operator onto $Y$ in the metric induced by $D$ is defined as $P_{Y,D}(x)=\underset{z\in Y}{\argmin} \ \|z-x\|^2_D$ for all $x\in\mathcal{H}$.

\medskip

We now formulate the general optimization problem we aim to solve and recall some preliminary technical results useful for the following convergence analysis.

We are interested in solving the problem:
\begin{equation}\label{minf}
{\min_{\x\in\H}} \ F(\x) \equiv f(\x) + g(\x) \, ,
\end{equation}
where
\begin{itemize}
\item $f:\H\rightarrow \R$ is convex and continuously differentiable with $L_f$~-Lipschitz continuous gradient on a closed convex set $Y\neq \emptyset$ with $\operatorname{dom}(g)\subseteq Y\subseteq \operatorname{dom}(f)$;
\item $g:\H\rightarrow \R\cup\{+\infty\}$ is proper,  convex, and lower semicontinuous.
\end{itemize}
Given $\bar{x}\in Y$, we also introduce the function $h_{\tau,D}(\cdot;\bar{x}):\H\rightarrow \R\cup\{+\infty\}$ defined for all $z\in\H$ as
\begin{equation}\label{eq:function_h}
h_{\tau,D}(z;\bar{x}):=f(\bar{x})+ \langle\nabla f(\bar{x}),z-\bar{x}\rangle+\frac{1}{2\tau}\|z-\bar{x}\|_D^2+g(z).
\end{equation}
Since $h_{\tau,D}(\cdot;\bar{x})$ is  strongly convex with respect to the $D$-norm with modulus $\frac{1}{\tau}$, it  has a unique minimizer $\hat{x}\in\mathcal{H}$, called the proximal--gradient point, which is given by
\begin{equation}\label{eq:prox-grad}
\hat{x}:=\prox_{\tau g}^D(\bar{x}-\tau D^{-1}\nabla f(\bar{x}))=\underset{z\in\H}{\operatorname{argmin}} \ h_{\tau,D}(z;\bar{x}).
\end{equation} 

In order to take into account possible inexact computations of $\hat{x}$, we now introduce ia suitable notion of approximation defined in terms of a fixed positive tolerance parameter (see \cite{Salzo-Villa-2012} for a detailed study).

\begin{definition}  \label{def:eps_approx}
Given $\bar{x}\in Y$, $\tau>0$, $D\in \mathcal{D}_{\eta_{inf}}^{\eta_{sup}}$ and $\epsilon\geq 0$, we say that a point $\tilde{x}\in\dom(g)$ is an $\epsilon-$approximation of the proximal--gradient point $\hat{x}$ and write $\tilde{x}\approx_\epsilon \hat{x}$ if 
\begin{equation}\label{eq:eps_approx}
h_{\tau,D}(\tilde{x};\bar{x})-h_{\tau,D}(\hat{x};\bar{x})\leq \epsilon.
\end{equation}
\end{definition}

The following scaled and inexact descent inequality holds (see also \cite[Lemma 2.3]{Bonettini2018a}).

\begin{Lemma}[Lemma 2.3, \cite{SAGE-FISTA}]\label{lem:technical}
Given $\bar{x}\in Y$, $\tau>0$, $D\in \mathcal{D}_{\eta_{inf}}^{\eta_{sup}}$, $\epsilon\geq 0$ and $\tilde{x}\approx_\epsilon \hat{x}$,  the following inequality holds for all $x\in\mathcal{H}$:
\begin{align}\label{eq:ine_fund}
F(\tilde{x})& +\frac{\|x-\tilde{x}\|_D^2}{2\tau}+\left(\frac{\|\tilde{x}-\bar{x}\|_D^2}{2\tau}-\mathbb{D}_f(\tilde{x},\bar{x})\right) \\ 
& \leq F(x)+\frac{\|x-\bar{x}\|_D^2}{2\tau}
+\epsilon+\frac{\sqrt{2\epsilon\tau}}{\tau}\|x-\tilde{x}\|_D.
\end{align}
\end{Lemma}
Descent inequalities in the form \eqref{eq:ine_fund} are the crucial tool in the analysis of convergence properties of proximal algorithms, as we will see in the following.

\section{Inexact S-FISTA: description and convergence result}

\subsection{Inexact S-FISTA}
We now describe the inexact, scaled and adaptive FISTA algorithm, dubbed inexact S-FISTA, for the solution of the composite convex optimization problem \eqref{minf} by means of a scaled and inexact inertial forward-backward splitting endowed with an adaptive backtracking strategy. Inexact S-FISTA is characterised by the following features:
\begin{itemize}
    \item the use of a variable metric in  \eqref{eq:prox-grad}, which is induced along the iterations by a sequence of linear, bounded and self-adjoint positive operators $\{D_k\}_{k\in\mathbb{N}}$, typically chosen so as to capture second order information of the differentiable part $f$ at the current iterate $x^{(k)}$ (see e.g. \cite{Bonettini-etal-2009,Bonettini2018a,Bonettini2019,Lanteri-etal-2001});  
    \item the inexact computation of the proximal-gradient point according to Definition \ref{def:eps_approx};
    \item a non-monotone backtracking strategy analogous to the one in \cite{Calatroni-Chambolle-2019}, which allows for possible increasing and decreasing of the step-size $\tau_{k+1}$ at each iteration; this strategy is  particularly helpful when the initial $\tau_0$ is chosen to be extremely small, which corresponds to a pessimistic estimate $L_0$ of $L_f$.
\end{itemize}
The proposed inexact S-FISTA is reported in Algorithm \ref{alg:S-FISTA}.

\begin{algorithm}[h]
\caption{Inexact S-FISTA($x^0$,$\tau_0$,$f$,$g$)}
\label{alg:S-FISTA}
\small{
\textbf{Parameters:} $\rho\in(0,1)$, $\delta\in(0,1]$, $\{\eta_{inf}^k\}_k$, $\{\eta_{sup}^k\}_k$ s.t. $0<\eta_{inf}\leq\eta_{inf}^k\leq \eta_{sup}^k\leq \eta_{sup}$.\\
\textbf{Initialization:} $x^{(-1)}=x^{(0)}$, $t_0\geq 1$,  $D_0\in\mathcal{D}_{\eta_{inf}^0}^{\eta_{sup}^0}$.

{\textsc{FOR $k=0,1,\ldots$ REPEAT}}
\begin{itemize}
    \item[] Choose $D_{k+1}\in\mathcal{D}_{\eta_{inf}^k}^{\eta_{sup}^k}$ and set $\tau_{k+1}^0=\frac{\tau_k}{\delta}$.\\
    {\textsc{FOR $i=0,1,\ldots $ REPEAT}}
    \begin{itemize}[leftmargin=1cm]

 \item[\textsc{1.}] $\tau_{k+1}=\rho^i \tau_{k+1}^0$
        \item[\textsc{2.}]  $t_{k+1}=\frac{1+\sqrt{1+4\frac{\tau_k}{\tau_{k+1}}t_k^2}}{2}$
        \item[\textsc{3.}] $y^{(k+1)}=\operatorname{proj}_{Y,D_{k+1}}\left(x^{(k)}+\left(\frac{t_k-1}{t_{k+1}}\right)(x^{(k)}-x^{(k-1)})\right)$
        \item[\textsc{4.}] Choose ${\epsilon}_{k+1}$ and compute $x^{(k+1)}\in\text{dom}(g)$ s.t. $x^{(k+1)}\approx_{\epsilon_{k+1}} \hat{x}^{k+1}$ with
\begin{equation*}
            \hat{x}^{(k+1)}:=\operatorname{prox}_{\tau_{k+1}g}^{D_{k+1}}(y^{(k+1)}-\tau_{k+1}D_{k+1}^{-1}\nabla f(y^{(k+1)}))
\end{equation*}
\normalsize
    \end{itemize}
    {\textsc{UNTIL}}
    $\mathbb{D}_f(x^{(k+1)},y^{(k+1)}) \leq\frac{1}{2\tau_{k+1}}\|x^{(k+1)}-x^{(k+1)}\|^2_{D_{k+1}}$
\end{itemize} 

{\textsc{UNTIL}} stopping criterion
}
\normalsize
\end{algorithm}

For $k\geq 0$ and as a first (preliminary) step, a linear, bounded and self-adjoint operator $D_{k+1}\in\mathcal{D}_{\eta_{inf}}^{\eta^k_{sup}}$ and a tentative step-size $\tau_{k+1}^0=\tau_k/\delta$ are chosen. Note that if $\delta<1$, then a larger step-size is attempted, similarly as in \cite{Calatroni-Chambolle-2019,Florea2020,Scheinberg-2014}, while if $\delta=1$ a classical Armijo-type backtracking is performed (see \cite{Beck-Teboulle-2009b}). Within the inner backtracking procedure indexed by $i= 0, 1,\ldots$, the quantities $t_{k+1}$ are updated depending on $\tau_{k+1}$. Then, the projected inertial point $y^{(k+1)}$ is computed via standard FISTA extrapolation. Finally, an approximated proximal-point $x^{(k+1)}\approx_{\epsilon_{k+1}} \prox_{\tau_{k+1} g}^{D_{k+1}}(y^{(k+1)}-\tau_{k+1}D_{k+1}^{-1}\nabla f(y^{(k+1)}))$ is computed in terms of the accuracy value $\epsilon_{k+1}\geq 0$ whose choice will be specified in Section \ref{sec:choice}. For each backtracking iteration, a check on the condition 
\begin{equation}\label{eq:backtracking}
\mathbb{D}_{f}(x^{(k+1)},y^{(k+1)})\leq \frac{\|x^{(k+1)}-y^{(k+1)}\|_{D_{k+1}}^2}{2\tau_{k+1}}
\end{equation}
is then performed . If \eqref{eq:backtracking} is not satisfied, then the step-size is reduced by a factor $\rho$ and new choices of $t_{k+1}, y^{(k+1)}$ are performed until \eqref{eq:backtracking} is satisfied. 

Note that when $\eta_{inf}^k=\eta_{sup}^k=1$, $D_k\equiv \mathcal{I}$ and $Y=\mathcal{H}$, Algorithm \ref{alg:S-FISTA} reduces to the adaptive variant of FISTA proposed in \cite{Scheinberg-2014}. 
When instead $D_k\neq \mathcal{I}$ and $\delta=1$, it corresponds to the inexact scaled forward-backward extrapolation method equipped with Armijo-type backtracking proposed in \cite{Bonettini2018a}. In this case, standard $\mathcal{O}(1/k^2)$ convergence rates can be proved, coherently with the result obtained in \cite[Theorem 3.1]{Bonettini2018a}.

Algorithm \ref{alg:S-FISTA} thus combines all the good features (variable scaling, adaptive backtracking and inexactness) of the algorithms in \cite{Beck-Teboulle-2009b,Scheinberg-2014,Bonettini2018a}, thus resulting in a general and flexible algorithm adapted to solve several imaging problems, see Section \ref{sec:results}.

\subsection{Convergence results}

We sketch in the following the main assumptions and convergence result proved in \cite[Section 3.1]{SAGE-FISTA} in the general strongly convex scenario and specified here for the convex optimization problem \eqref{minf}. We start specifying a technical assumption on the sequence of operators $\{D_k\}_{k\in\mathbb{N}}$ which has been previously employed in analogous works to prove the convergence of variable metric FB algorithms, see, e.g., \cite{Combettes-Vu-2014,Bonettini-Loris-Porta-Prato-2015,Bonettini-Porta-Ruggiero-2016,Bonettini2018a}.

\medskip

\begin{assumption}\label{ass:1}
There exists a sequence of real nonnegative numbers $\{\gamma_k\}_{k\in\mathbb{N}}$ s.t. $\sum_{k=0}^\infty\gamma_k<+\infty$ and, for all $k\geq 0$, the following condition holds  
\begin{eqnarray} \label{eq:assumption}
D_{k+1}&\preceq&(1+\gamma_{k+1})D_k,\label{eq:Dk_cond}.
\end{eqnarray}
\end{assumption}

Note, in particular, that proceeding in \cite[Remark 3.3]{SAGE-FISTA} one can show that Assumption \ref{ass:1} holds when $D_k\in\mathcal{D}_{\eta_{inf}^k}^{\eta_{sup}^k}$ for all $k\geq 0$, and
\begin{equation}\label{eq:eta_def}
\eta_{inf}^k=\eta-\nu_{inf}^k, \quad \eta_{sup}^k=\eta+\nu_{sup}^k,
\end{equation}
where $0\leq\nu_{inf}^k< \eta$, $\nu_{sup}^k\geq 0$, and $\sum_{k=0}^{\infty}\nu_{inf}^k<+\infty$, $\sum_{k=0}^{\infty}\nu_{sup}^k<+\infty$, i.e. the upper and lower bounds converge to the same positive value at a sufficiently fast rate. 
In applications where $\mathcal{H}=\mathbb{R}^n$, $\|\cdot\|=\|\cdot\|_2$, and $\{D_k\}_{k\in\mathbb{N}}$ are diagonal matrices, it is always possible to impose condition \eqref{eq:eta_def} by forcing the diagonal elements of $D_k$ to belong to the interval $[\eta-\nu_{inf}^k,\eta+\nu_{sup}^k]$ (see Section \ref{sec:choice}). By doing so, the scaling matrices tend to $\eta\mathcal{I}$ as iterations progress.

\medskip

We now state the main convergence result for Algorithm \ref{alg:S-FISTA}. Its proof is based on the use of the descent inequality \ref{lem:technical} and of technical results based on induction arguments. We refer the reader to \cite[Theorem 3.1]{SAGE-FISTA} for the general statement and proof of this result in the possibly strongly convex setting. For more details on the convex case, we also refer to \cite[Corollary 3.3]{SAGE-FISTA}.

\medskip

\begin{theorem}[Theorem 3.1 \& Corollary 3.3 \cite{SAGE-FISTA}]\label{th:convergence}
Let $x^*\in\H$ be a solution of \eqref{minf}. Suppose that Assumption \ref{ass:1} holds and that the sequence $\{\epsilon_k\}_{k\in\mathbb{N}}$ is chosen as
\begin{equation}\label{eq:prefixed3}
\epsilon_{k+1}=\begin{cases}
\mathcal{O}(a^{k+1}), \quad &\text{if } \ \delta <1\\
\displaystyle \frac{b_{k+1}}{(k+1+t_0)^2}, \quad &\text{if } \ \delta =1
\end{cases},
\end{equation}
where $a<\delta$ and $\{b_k\}_{k\in\mathbb{N}}$ is such that $b_k\geq 0$ and $b:=\sum_{k=0}^{\infty}\sqrt{b_k}<\infty$. Let further define  $\gamma:=\lim_{k\to\infty}\prod_{i=1}^k (1+\gamma_i)<+\infty$. Then, for all $k\geq 0$, we have
\begin{equation*}
F(x^{(k+1)})-F(x^*)\leq \frac{C}{(k+1+t_0)^2},
\end{equation*} 
where $C=C(x^*,x^{(0)},\tau_0,D_0,t_0,\rho,\delta,L_f,\eta_{inf},\eta_{sup},\gamma, a, b)$ is a constant depending on the algorithmic parameters.
\end{theorem}

\medskip

The $\mathcal{O}(1/k^2)$ convergence rate result for the function values is well-known for FISTA \cite{Beck-Teboulle-2009b} and for its adaptive \cite{Scheinberg-2014} and inexact and variable metric variants \cite{Schmidt2011,Bonettini2018a}. Theorem \ref{th:convergence} unifies the existing results  under one general result. Note that as noted in \cite{Calatroni-Chambolle-2019} and in \cite[Remark 3.4]{SAGE-FISTA}, the same rate of convergence can be obtained by avoiding the dependence of $C$ on the possibly unknown constant $L_f$ in terms of the average quantity:
\begin{equation}   \label{eq:Lbar}
\sqrt{\bar{L}_{k+1}} := \frac{1}{\frac{1}{k+2} \sum\limits_{i=0}^{k+1}  \sqrt{\tau_i}},
\end{equation}
which requires the storage of all values $\tau_i$ computed along the iterations.

\section{Inexact proximal points and variable metric selection}  \label{sec:choice}

We shortly describe in this section how the inexactness condition \ref{eq:eps_approx} and the assumption \eqref{eq:eta_def} (guaranteeing \eqref{eq:assumption}) for the variable scaling operators can be implemented in practice.

\subsection{Inexact computation of proximal points with $\epsilon_k$-accuracy}  \label{sec:inexact}

We recall here the general strategy detailed in \cite[Section 4.2]{Bonettini2018a} for computing an inexact proximal--gradient point guaranteeing the condition \eqref{eq:eps_approx} in Definition \ref{def:eps_approx}. For the following examples, we will require that the  function $g$ in \eqref{minf} can be expressed in the form:

\begin{equation}  \label{eq:g_form}
g(x) = \sum_{i=1}^p \phi_i (M_ix) + \psi(x),
\end{equation}
where $M_i:\mathcal{H}\to \mathcal{Z}_i$ are linear bounded operators between Hilbert spaces and $\phi_i:\mathcal{Z}_i\to \mathbb{R}\cup\left\{+\infty\right\}$, $\psi:\mathcal{H}\to \mathbb{R}\cup\left\{+\infty\right\}$ are proper,  convex and lower semicontinous functions.
For image restoration problems, the couple $(\phi_i, M_i)$ is typically associated to the use of gradient-type regularization terms (e.g., by setting $\phi_i:\mathbb{R}^2\to\mathbb{R}, \phi_i(v) = \lambda\|v\|$ and $M_i = \nabla_i$, the discrete image gradient), while the function $\psi$ may encode further requirements on the desired signal such as a positivity constraint or a convex perturbation term.
At each iteration $k\geq 1$ of inexact S-FISTA, the primal subproblem to be solved to compute $x^{(k)}$ takes the form:

\small{
\begin{equation}  \label{eq:primal_pb}
   \min_{x\in \mathcal{H}}~ \left\{ \mathcal{P}_{\tau_{k},D_{k}}(x) := \sum_{i=1}^p \phi_i (M_ix) + \psi(x) + \frac{1}{2\tau_{k}} \|x- \bar{y}^{(k)} \|^2_{D_{k}} \right\},
\end{equation}}
\normalsize
where $\bar{y}^{(k)}:=y^{(k)}-\tau_{k} D_{k}^{-1}\nabla f (y^{(k)})$. Under suitable assumptions (see, e.g., \cite[Section 4.1]{Bonettini2018a} and reference therein), solving \eqref{eq:primal_pb} is equivalent to maximizing the associated dual function $ \mathcal{Q}_{\tau_{k},D_{k}}:\mathcal{Z}_1\times\ldots\times \mathcal{Z}_p\rightarrow \mathbb{R}\cup\{+\infty\}$ obtained by Fenchel conjugation of the functions $\phi_i$, $i=1,\ldots,p$. Moreover, one can deduce the following inequality
\begin{align*}
& h_{\tau_k,D_k}(x;\bar{y}^{(k)})-h_{\tau_k,D_k}(\prox_{\tau_{k}\psi}^{D_{k}} (\bar{y}^{(k)});\bar{y}^{(k)}) \\
& \leq   \mathcal{P}_{\tau_{k},D_{k}}(x) - \mathcal{Q}_{\tau_{k},D_{k}}(w), 
\end{align*}
for all $w\in\mathcal{Z}$, which entails that a sufficient condition for a point $x$ to be an $\epsilon_{k}$-approximation as in Definition \ref{def:eps_approx} is the existence of a dual point $w$ such that $\mathcal{P}_{\tau_{k},D_{k}}(x) - \mathcal{Q}_{\tau_{k},D_{k}}(w)\leq \epsilon_{k}$. 
Assuming that $g$ is continuous on $\text{dom}(g)=\text{dom}(\psi)$, an $\epsilon_{k}$-approximation can thus be computed by defining a dual sequence $\left\{ w ^{(k,l)}\right\}_{l\in\mathbb{N}}\subset \mathcal{Z}$ converging to the solution of the dual problem and a corresponding primal sequence $\left\{ x ^{(k,l)}\right\}_{l\in\mathbb{N}}\subset \mathcal{H}$ defined for $l\in\mathbb{N}$ by
\begin{equation}
        x^{(k,l)} := \prox_{\tau_{k}\psi}^{D_k}\left(\bar{y}^{(k)} - \tau_{k}D_k^{-1} M^* w^{(k,l)} \right),\label{eq:primal_seq}
\end{equation}
and then stopping the iterates whenever
\begin{equation}  \label{eq:criterion_eps}
\mathcal{P}_{\tau_{k},D_k}(x^{(k,l)}) - \mathcal{Q}_{\tau_{k},D_k}(w^{(k,l)})\leq \epsilon_k.
\end{equation}
Following \cite[Proposition 4.2]{Bonettini2018a} we have that the procedure is well-defined.
The dual sequence $\left\{w^{(k,l)} \right\}_{l\in\mathbb{N}}$ can be generated using an efficient inner FISTA routine, provided that the extrapolation parameters are chosen in a way that weak convergence of the iterates is guaranteed (see, e.g.~ \cite{Chambolle-Dossal-2014}), whereas the primal sequence \eqref{eq:primal_seq} can be computed in closed form in many practical situations, such as the ones reported in the following sections.

\subsection{Split-gradient strategy for variable metric selection}

For choosing the scaling matrices $\left\{D_k\right\}_{k\in\mathbb{N}}$, we exploit the split-gradient strategy proposed in \cite{Lanteri-etal-2001} and later used in several works (see, e.g., \cite{Bonettini-Loris-Porta-Prato-2015,Bonettini2018a}) for which the decomposition $-\nabla f(x) = U(x) - V(x)$ with $U(x)\geq 0$ and $V(x)>0$ and the choice $D_k=\text{diag}\left( y^{(k)}/V(y^{(k)}) \right)^{-1}$ is made.  
In order to ensure the conditions required by Assumption \ref{ass:1}, we further need to introduce a thresholding parameter $\gamma_k$, thus considering
\begin{equation}  \label{eq:scaling_general}
    D_k = \text{diag}\left(\max\left(\frac{1}{\gamma_k},\min\left(\gamma_k,  \frac{y^{(k)}}{V(y^{(k)})} \right)\right)\right)^{-1},
\end{equation}
for thresholding parameters $\gamma_k$ defined by:
\begin{equation} \label{eq:gamma_k_pb1}
    \gamma_k = \sqrt{1+\frac{s_1}{(k+1)^{s_2}}},\quad \text{where }s_1 >0, \ s_2>1.
\end{equation}
As it is obvious, the choice \eqref{eq:scaling_general} depends on the specific problem considered, due to presence of the the function $V(\cdot)$. Note that when $s_1=0$, then $D_k\equiv\mathcal{I}$, so the standard Euclidean metric is recovered. Numerically, it is good practice to choose a large value of $s_1$ to benefit from the use of the ``Newton-type"  metrics in the early iterations of the algorithm, while letting $s_2$ drive the asymptotic behaviour.

\section{Numerical results} \label{sec:results}

We now apply the Inexact S-FISTA Algorithm \ref{alg:S-FISTA} to solve an examplar image deblurring problem. As observed in \cite{Bertero2018} , the use of variable metric optimisation algorithms has been showed to render particularly effective in the case of data corrupted by signal-dependent Poisson noise, which is frequently encountered in microscopy and astronomical imaging.  In our experiments Poisson noise is simulate by means of the MATLAB \texttt{imnoise} routine.

\subsection{Problem formulation}

For a given image $z\in \mathbb{R}_{\geq 0}^n$, we thus consider the ill-posed image restoration problem 
$$
\text{find}\quad x\in\mathbb{R}_{\geq 0}^n\quad\text{s.t.}\quad z = \mathcal{P}\left(Hx + b \right),
$$
where $H\in\mathbb{R}^{n\times n}$ is the blurring operator computed for a given Gaussian Point Spread Function (PSF) with standard deviation $\sigma_{PSF}>0$, the term $b\in\mathbb{R}^n_{> 0}$ is a positive background term and $\mathcal{P}(w)$ models Poisson noise degradation. As it is well-known by standard Maximum A Posteriori (MAP) estimation, the data fidelity term modelling the presence of Poisson noise is the generalized Kullback-Leibler (KL) divergence functional defined by:
\small{
\begin{equation}  \label{eq:KL_def}
    KL(Hx+b;z) := \sum_{i=1}^n \left( z_i \log \frac{z_i}{(Hx)_i + b_i} + (Hx)_i + b_i - z_i \right),
\end{equation}}
\normalsize
where the convention $0 \log 0 = 0$ is adopted. Note that if the operator $H$ has nonnegative entries and if it has at least one strictly positive entry for each row and column (i.e. 
$He>0, H^Te>0$ for $e\in\mathbb{R}^n$ being the vector of all ones), the function $KL(H\cdot+b;z)$ is nonnegative, convex and coercive on the non-negative orthant $Y:= \left\{ x\geq 0 \right\}$, see, e.g., \cite{Harmany12}. 
We couple \eqref{eq:KL_def} with the non-smooth isotropic Total Variation (TV) semi-norm which is often employed for solving imaging problems due to its edge-preserving properties:
\begin{equation}   \label{eq:TV}
TV(x) := \| \nabla x \|_{2,1} = \sum_{i=1}^n \| \nabla_i x \|_2,
\end{equation}
where, for $i=1,\ldots,n$ we denote by $\nabla_i$  the standard forward-difference  image gradient operator. 
Reflexive boundary conditions for the computation of such discretization are used so that for $H$  and $H^T$ matrix-vector products can be efficiently computed in terms of the Discrete Cosine Transform (DCT). By further imposing a non-negativity constraint on the orthant $Y$, we thus  end up with the following composite optimization problem:
\begin{equation}  \label{pb:TV_KL}
    \min_{x\in\mathbb{R}^n}~ KL(Hx+b;z) + \lambda TV(x) + \iota_Y(x)
\end{equation}
where $\lambda>0$ is a regularization parameter and where $\iota_Y(\cdot)$ stands for the indicator function of  $Y$.

\begin{figure}[t!]
\begin{center}
\begin{tabular}{ccc}
    \includegraphics[scale=0.15]{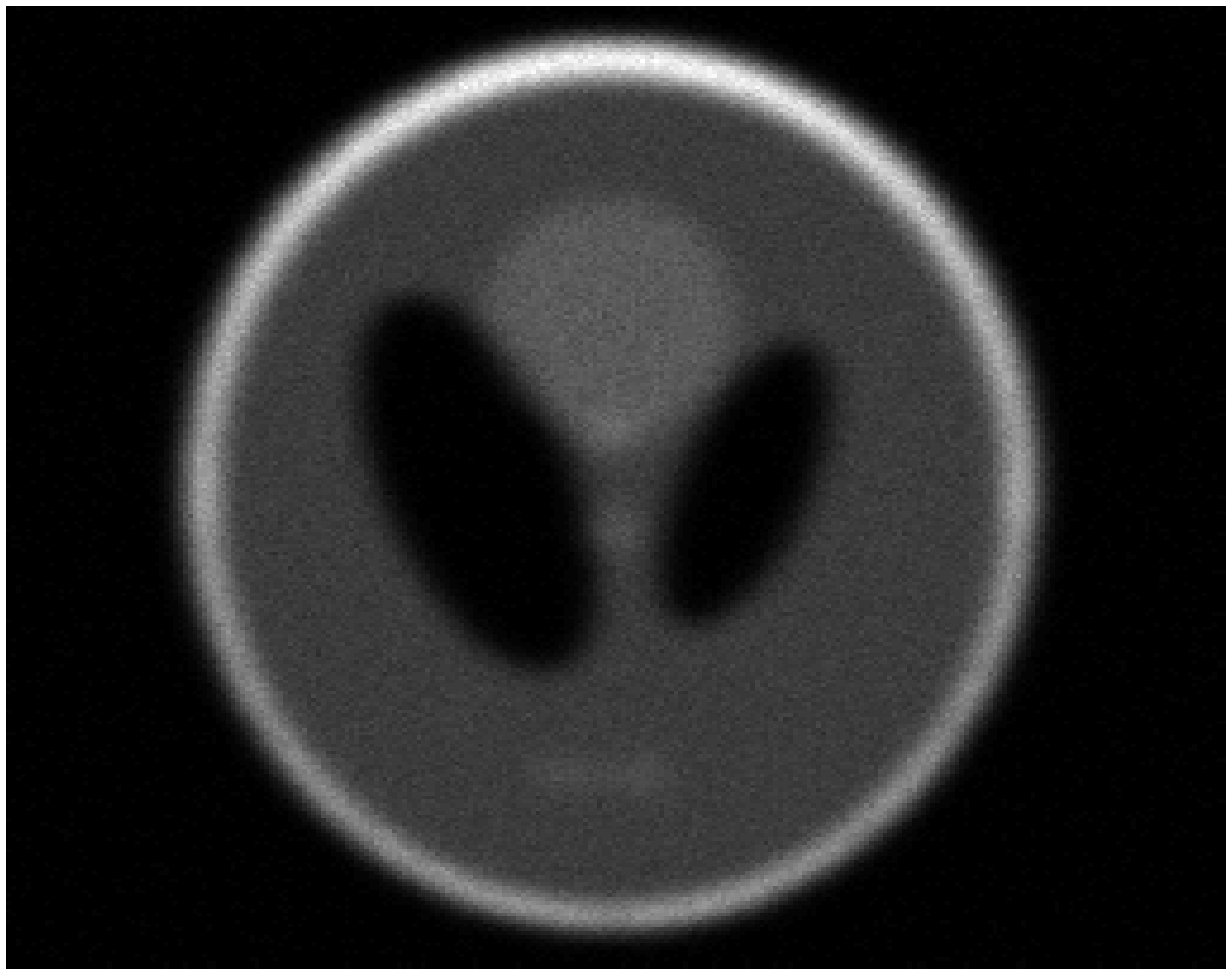} 
  &
    \includegraphics[scale=0.15]{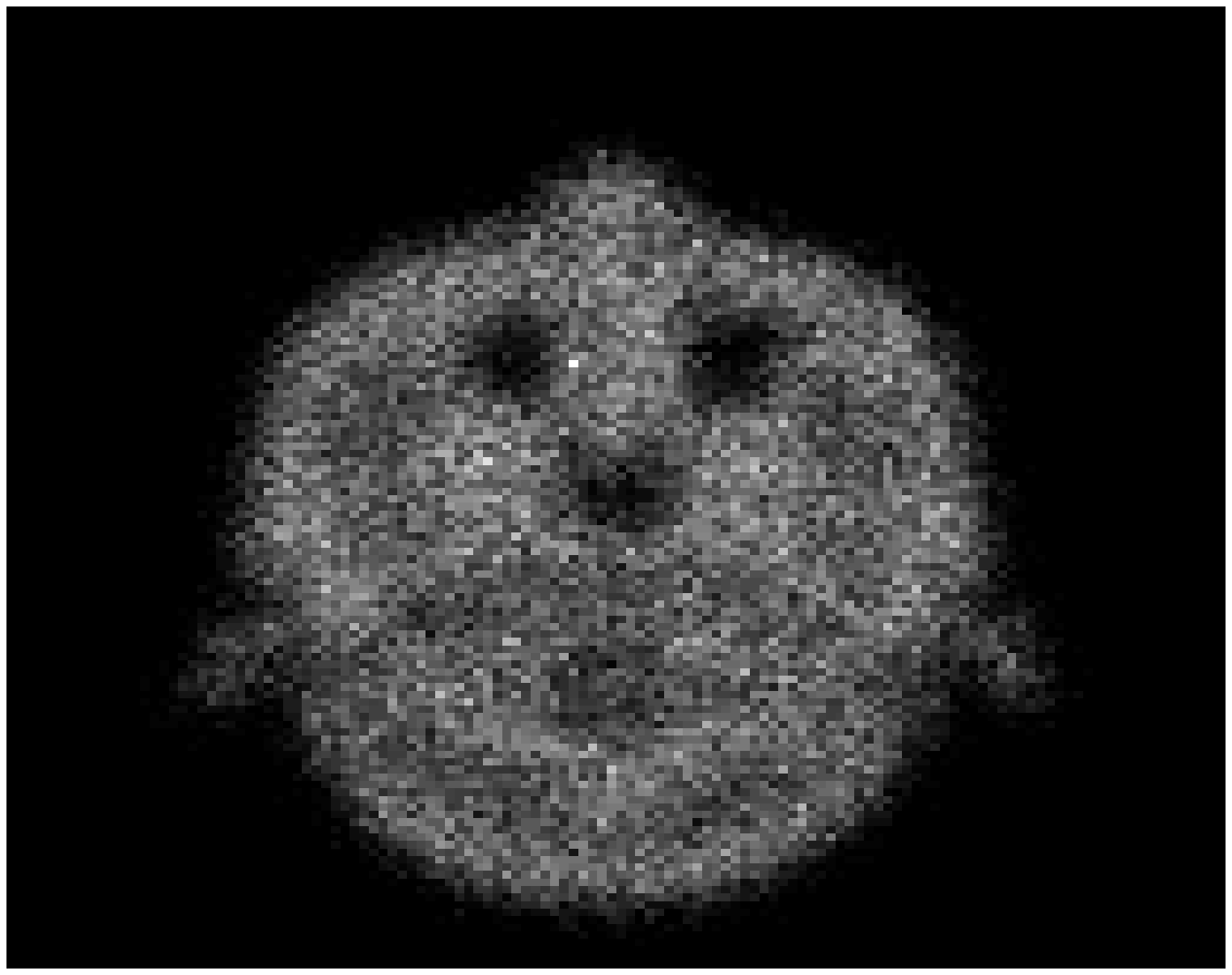} \\
  (a) \texttt{phantom}
 &
    (c) \texttt{mri}
\end{tabular}
    \caption{Blurred and noisy images used in the tests.}
    \label{fig:test_images}
\end{center}
\end{figure}

\begin{figure}[t!]
\begin{center}
		\scriptsize{
			\begin{tabular}{|c|c|c|c|c|c|c|}
				\hline
				\textbf{Image}   & \textbf{Size}   & \textbf{Range} & \textbf{$\sigma_{PSF}$} & \textbf{$b$} & \textbf{$\lambda$} &  \textbf{$L_f$} \\ \hline\hline
				\texttt{phantom} & $256\times 256$ & $[1,878]$        & $1.4$                   & $10$    & $0.004$            &  $8.78$       \\ \hline
				\texttt{mri}     & $128\times 128$ & $[0,170]$      & $3.2$                   & $0.5$        & $0.015$            &  $680$     \\ \hline
			\end{tabular}%
			\captionof{table}{Model parameters  for the test images in Figure.}
			\label{table:details_figs}
		}
\end{center}
\end{figure}

Recalling \eqref{minf}, we thus fix:
\[
f(x)=KL(Hx+b;z),\qquad g(x)=\lambda TV(x)+\iota_{Y}(x)
\]
so that:
\begin{align}  
    \nabla f(x) & = H^T e - H^T \left( \frac{z}{Hx + b}\right), \label{eq:parameters_PoisDeb}\\ L_f& =\frac{\max z_i}{b^2} \max (H^T e) \max(H e),\nonumber
\end{align}
where $e\in\mathbb{R}^n$ is a vector of ones. The expression of $\nabla f$ shows that the estimation of $L_f$  depends on the range on the data $z$, which in our examples is thus intentionally allowed to vary, i.e. images are not normalised within a fixed range. We report in Figure \ref{fig:test_images} and Table \ref{table:details_figs} the figures and the numerical details of the blurred and noisy images used in our numerical experiments.

We pre-compute an approximation $x^*$ of the desired solution by running standard FISTA for 3000 iterations. To assess convergence, we compute the relative objective error $(F(x^{(k)}) - F(x^*))/F(x^*)$, $k\geq 1$ both along the iterations and with respect to (at most) the first 30 seconds of run. The backtracking parameters are set $\rho=0.85$ and $\delta\in\left\{1,0.98 \right\}$, depending on whether a classical Armijo-type or adaptive backtracking is performed. We further set a maximum of \texttt{maxiter}$=200$ outer and \texttt{max\_bt}$=10$ inner backtracking iterations and initialise $t_0=1$ and $x^{(0)}=z$. The initial values for $\tau_0$ are specified in the captions of the following results.

Regarding the sequence of scaling matrices $\left\{D_k\right\}_{k\in\mathbb{N}}$, we consider the diagonal split-gradient strategy \eqref{eq:scaling_general} which in this case corresponds to:
\begin{equation}  \label{eq:metric_KL}
    D_k = \text{diag}\left(\max\left(\frac{1}{\gamma_k},\min\left(\gamma_k, \frac{y^{(k)}}{H^T e}\right)\right)\right)^{-1},
\end{equation}
where, notice, that the matrix depends explicitly on the extrapolated point $y^{(k)}$ and
where the thresholding parameters $\gamma_k$ are defined as in \eqref{eq:gamma_k_pb1}. 

As far as the choice of the sequence $\left\{\epsilon_{k+1}\right\}_{k\in\mathbb{N}}$,
we observe that condition  \eqref{eq:prefixed3} in Theorem \ref{th:convergence} is guaranteed by choosing $a=\frac{\delta}{2}$ and $b_k=1/k^{2.1}$.  Under these choices, the inexact computation of the proximal operator of $g$ can thus be performed proceeding as in \ref{sec:inexact}, after noticing that the function $g(\cdot)$ can be cast in the form \eqref{eq:g_form} by choosing $\phi_i(v) = \lambda \|v\|_2, M_i =  \nabla_i, i=1,\ldots,n$ and $\psi(x) = \iota_{Y}(x) $, so that we have for all $ l\in\mathbb{N}$
\begin{align}
x^{(k,l)}  
 = P_{Y,D_k}  \left( y^{(k)}-\tau_k D_{k}^{-1}\left(\nabla f(y^{(k)})   + M^* w^{(k,l)} \right)  \right). \label{eq:primal_scaled}
\end{align}

\subsection{Numerical experiments}

\begin{figure}[t!]
\begin{tabular}{c@{}c}
    \includegraphics[height=3.3cm]{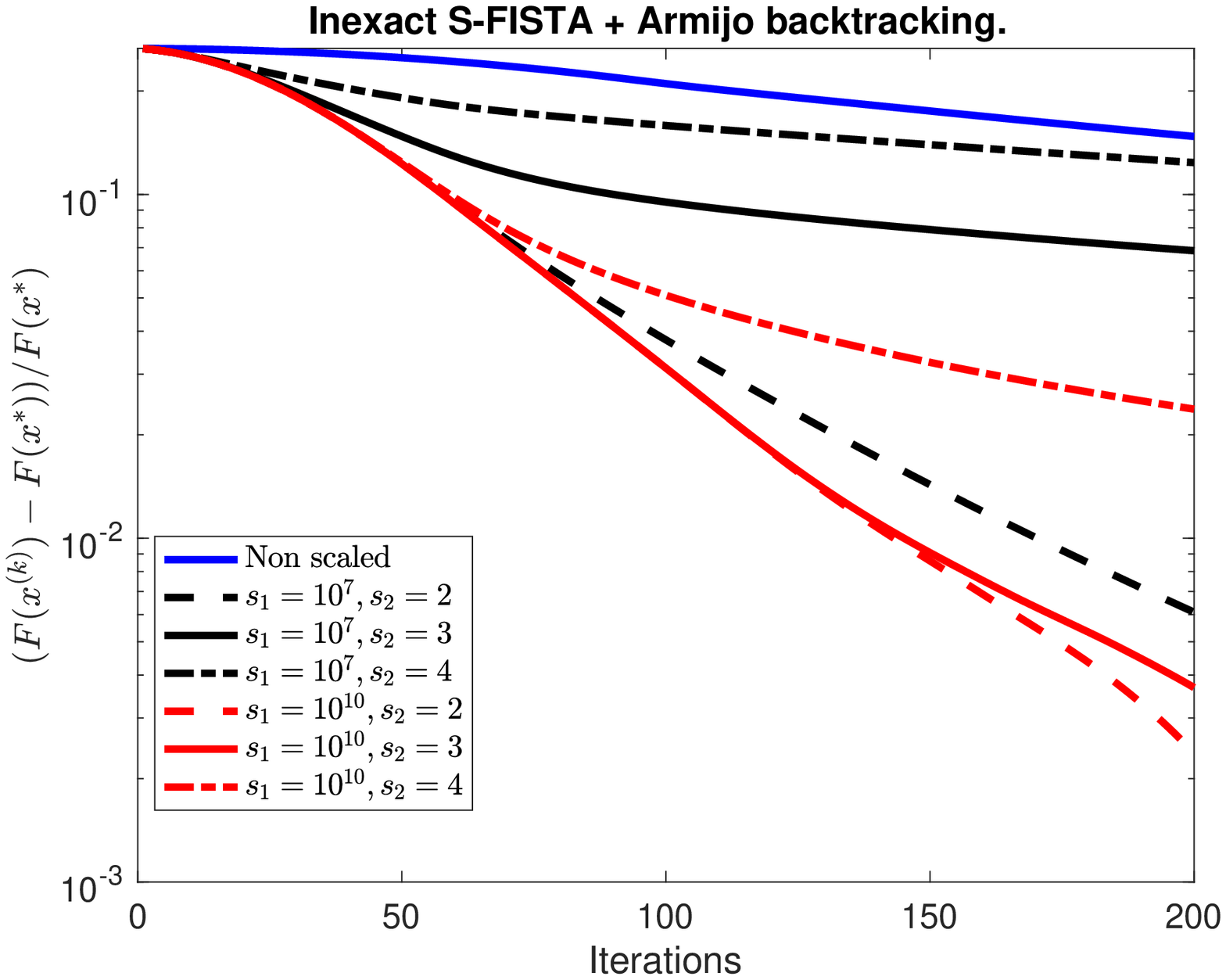} &  
    \includegraphics[height=3.3cm]{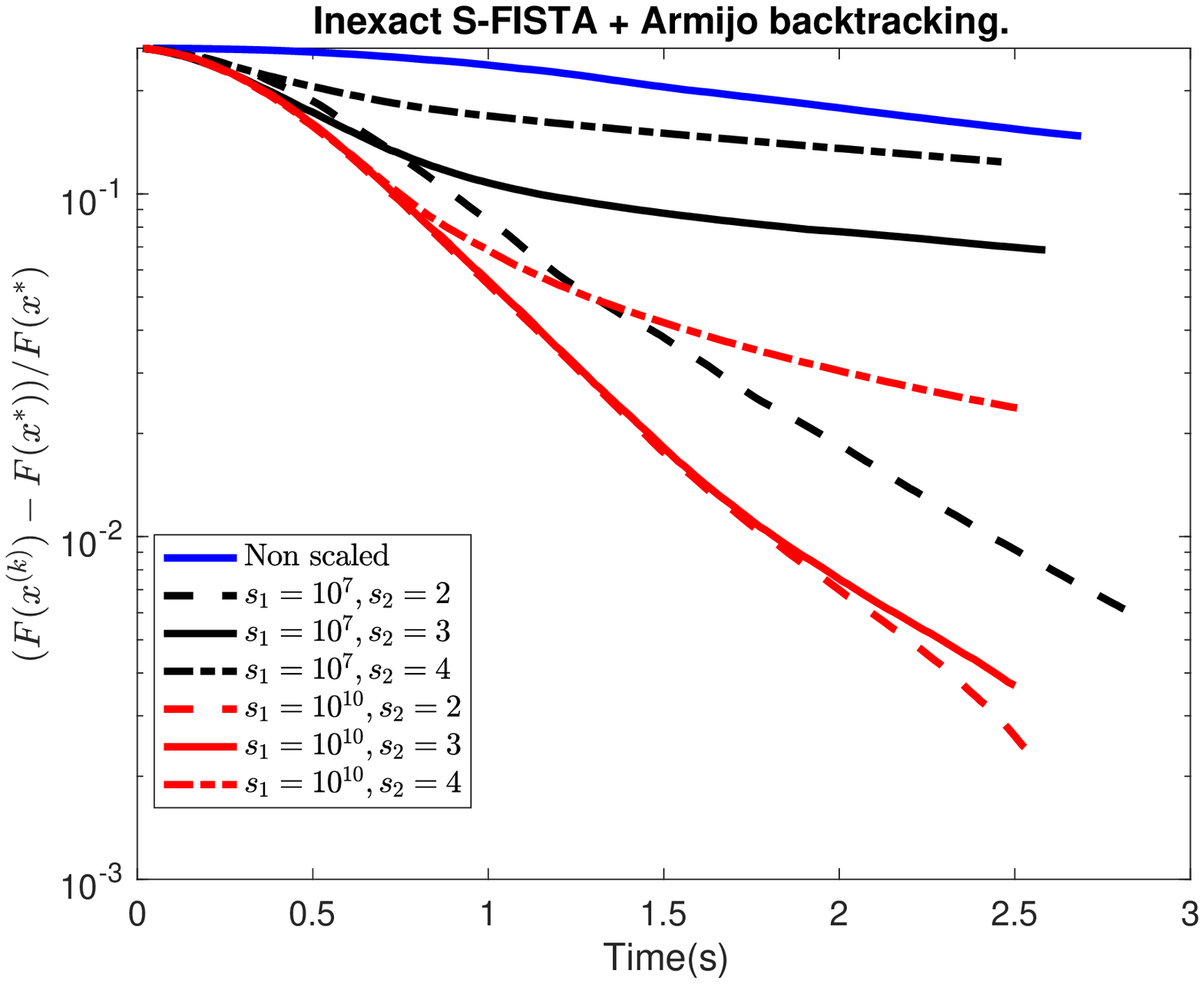} \\
     \includegraphics[height=3.3cm]{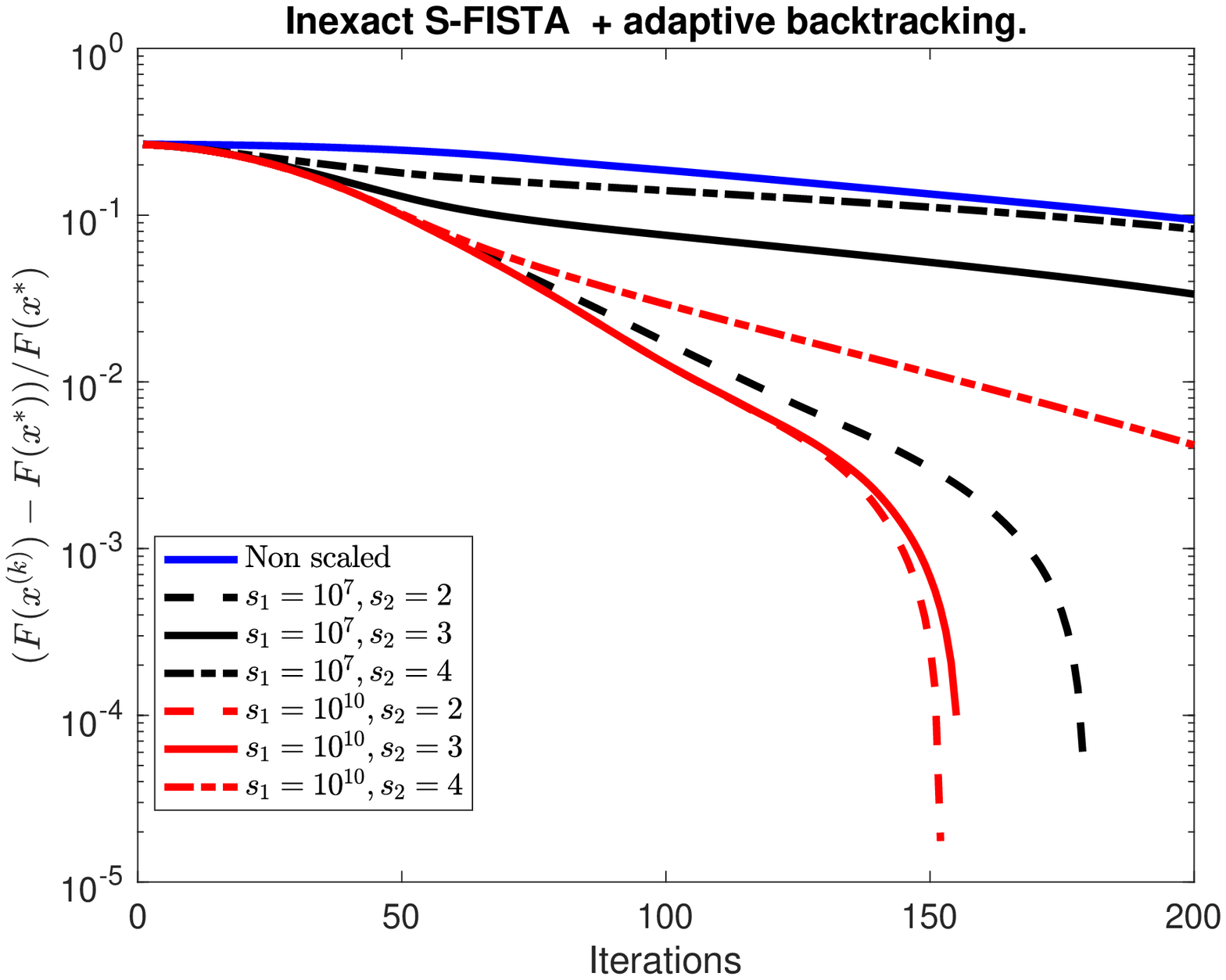} & \includegraphics[height=3.3cm]{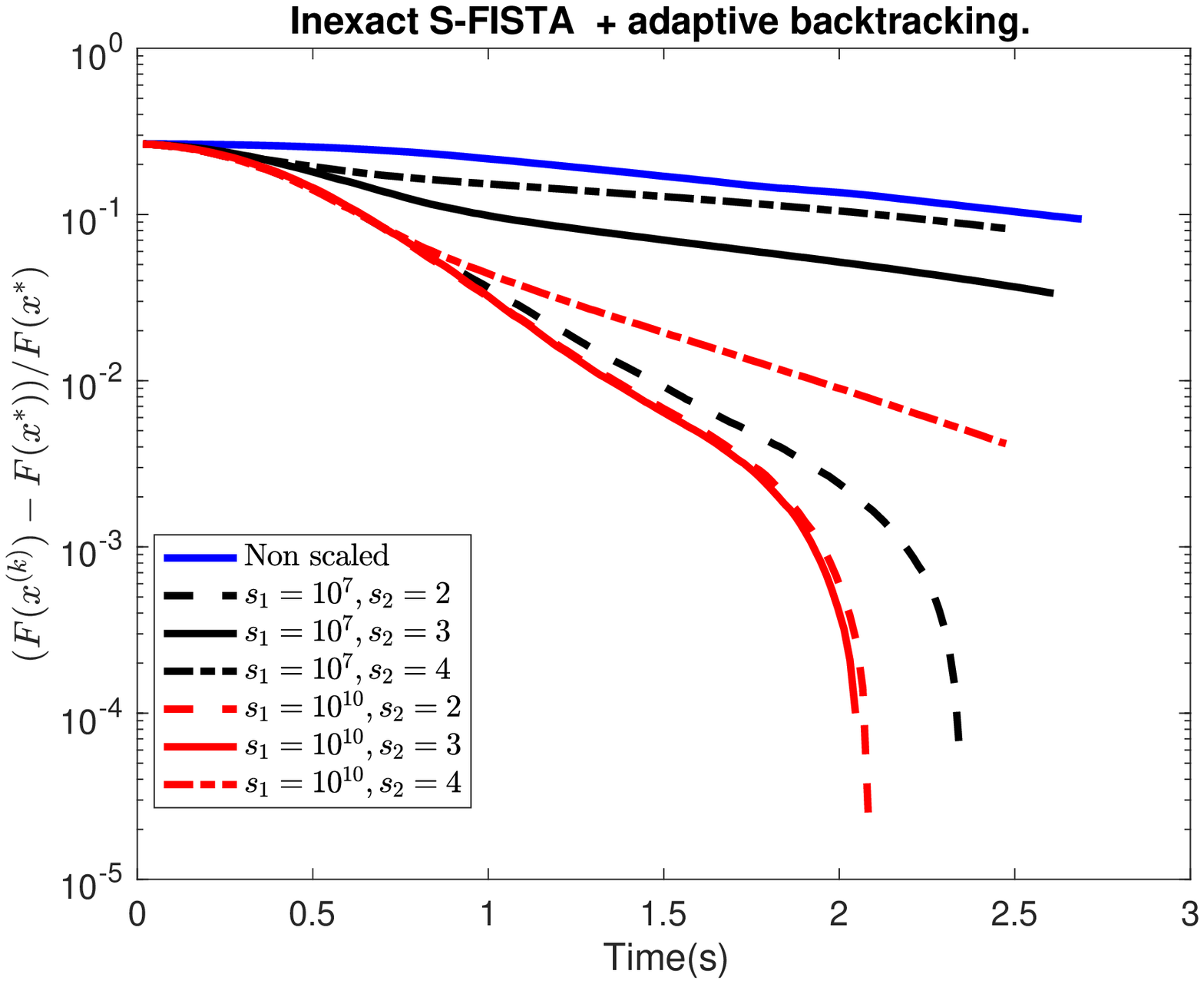} \\
     (a) Relative rates VS it. & (b) Relative rates VS CPU times. 
\end{tabular}
\caption{Inexact S-FISTA with Armijo (top row) and adaptive (bottom row) backtracking for different choices of $s_1$ and $s_2$ for problem \eqref{pb:TV_KL} on \texttt{mri} image, $L_0=200.$}
    \label{fig:mri}
\end{figure}

In Figure \ref{fig:mri} we report the results obtained by applying the Inexact S-FISTA Algorithm \ref{alg:S-FISTA} with Armijo and adaptive backtracking for different choices of the scaling parameters $s_1$ and $s_2$ in \eqref{eq:gamma_k_pb1} in correspondence of the \texttt{mri} test image for initial $L_0=200$. We observe that both in terms of convergence speed and computational times scaled algorithms outperform their non-scaled counterparts.

\begin{figure}[t!]
\begin{tabular}{c@{}c}
    \includegraphics[height=3.3cm]{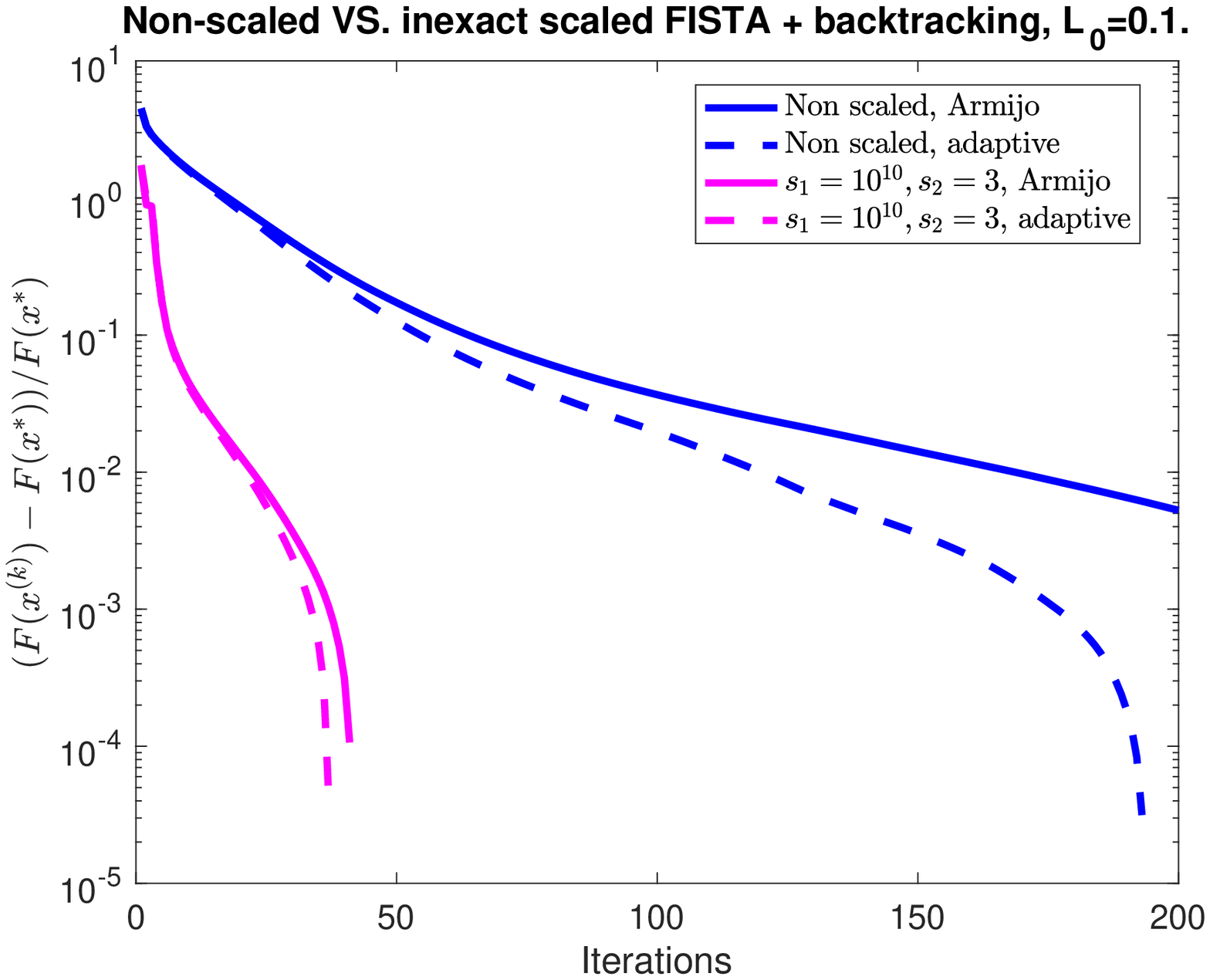}&  
    \includegraphics[height=3.3cm]{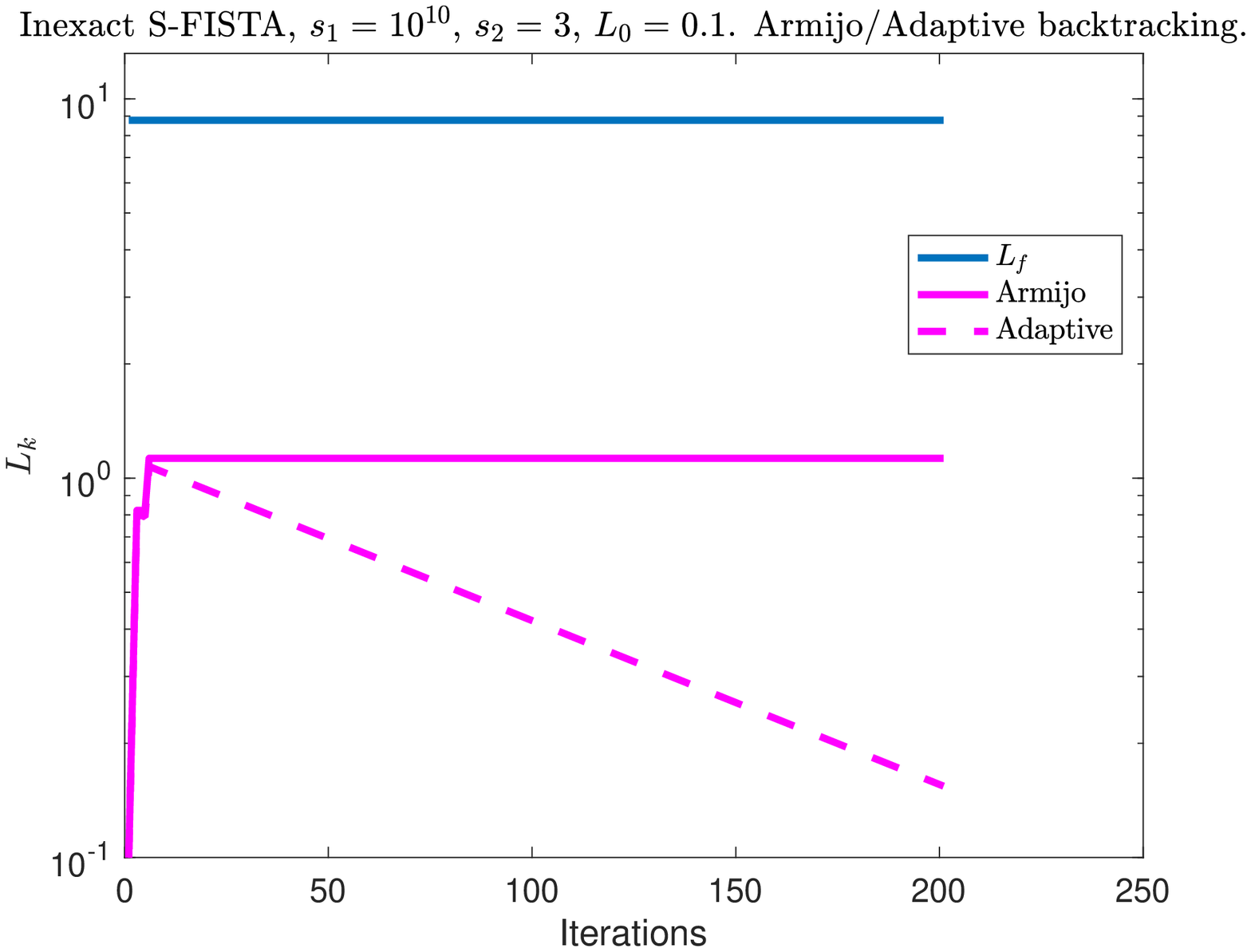} \\
     (a) Relative rates VS it. & (b) $L_k$ estimates 
\end{tabular}
\caption{(a): Relative convergence rates of non-scaled FISTA and Inexact S-FISTA with Armijo and adaptive backtracking for $s_1=10^{10}$ and $s_2=3$ for problem \eqref{pb:TV_KL} on \texttt{phantom} image, $L_0=0.01.$ (b) Lipschitz constant estimates $L_k$ for Inexact S-FISTA with Armijo/adaptive backtracking.}
    \label{fig:phantom}
\end{figure}

We run similar numerical tests on the \texttt{phantom} image with initial Lipschitz constant estimate $L_0=0.1$. In Figure \ref{fig:phantom} (a) we compare the rate of convergence along iterations of non-scaled FISTA algorithm and the Inexact S-FISTA algorithm for suitable choice of scaling parameters $s_1$ and $s_2$, both endowed with Armijo and adaptive backtracking. The combination of the variable scaling with the adaptive backtracking significantly improves convergence speed, thus allowing better precision. 
As far as the estimation of the Lipschitz constant $L_f$ is concerned, we finally report in Figure \ref{fig:phantom} (b) a comparison between the values $L_k$ estimated along the iterations via Armijo ($\delta=1$) and adaptive ($\delta<1$) backtracking. While performing analogously during the early iterations, due to its non-monotonicity, the adaptive backtracking strategy allows for local adjustments of the estimation which contributes to obtain faster convergence. Note, that this features allows also to correct possible too large misspecifications of $L_0$ which cannot be corrected by standard Armijo backtracking.

\section{Conclusions}

In this short note, we considered the convex instance of the SAGE-FISTA forward-backward algorithm recently proposed in \cite{SAGE-FISTA} to solve strongly convex composite optimisation problems. The proposed algorithm \ref{alg:S-FISTA} takes explicitly into account possible inexact evaluations of the proximal steps by means of a suitable notion of inexactness (Definition \ref{def:eps_approx}) defined in terms of a sequence of precision parameters $\left\{ \epsilon_k\right\}_{k\in\mathbb{N}}$. Moreover, it is defined in terms of a sequence of variable diagonal positive definite operators $\left\{ D_k\right\}_{k\in\mathbb{N}}$ satisfying a suitable ordering (Assumption \ref{ass:1}) and endowed of an adaptive backtracking strategy which allows for local adjustments of the algorithmic step-size. Theorem \ref{th:convergence} guarantees that for suitable choices of the sequence $\left\{ \epsilon_k\right\}_{k\in\mathbb{N}}$, the quadratic convergence rate for the function values, typically proved in FISTA-type algorithms, is still guaranteed. After describing in Section \ref{sec:choice} the practical implementation of the inexactness and variable metric steps of the algorithm, we report in Section \ref{sec:results} some numerical results performed on TV image deblurring problems for images corrupted with Poisson noise. Our experience shows that the combination of variable scaling and adaptive backtracking significantly improves convergence speed and favours better accuracy.

\bibliographystyle{abbrv}
\bibliography{biblio.bib}

\begin{thebibliography}{10}

\bibitem{Beck-Teboulle-2009b}
A.~Beck and M.~Teboulle.
\newblock A fast iterative shrinkage-thresholding algorithm for linear inverse
  problems.
\newblock {\em SIAM J. Imaging Sci.}, 2:183--202, 2009.

\bibitem{Bertero2018}
M.~Bertero, P.~Boccacci, and V.~Ruggiero.
\newblock {\em Inverse Imaging with {P}oisson {D}ata}.
\newblock 2053-2563. IOP, 2018.

\bibitem{Bonettini-Loris-Porta-Prato-2015}
S.~Bonettini, I.~Loris, F.~Porta, and M.~Prato.
\newblock Variable metric inexact line--search based methods for nonsmooth
  optimization.
\newblock {\em SIAM J. Optim.}, 26(2):891--921, 2016.

\bibitem{Bonettini2019}
S.~Bonettini, F.~Porta, M.~Prato, S.~Rebegoldi, V.~Ruggiero, and L.~Zanni.
\newblock {\em Recent Advances in Variable Metric First-Order Methods}, pages
  1--31.
\newblock Springer International Publishing, Cham, 2019.

\bibitem{Bonettini-Porta-Ruggiero-2016}
S.~Bonettini, F.~Porta, and V.~Ruggiero.
\newblock A variable metric forward-backward method with extrapolation.
\newblock {\em SIAM J. Sci. Comput.}, 38:A2558--A2584, 2016.

\bibitem{Bonettini2018a}
S.~Bonettini, S.~Rebegoldi, and V.~Ruggiero.
\newblock Inertial variable metric techniques for the inexact forward--backward
  algorithm.
\newblock {\em SIAM J. Sci. Comput.}, 40(5):A3180--A3210, 2018.

\bibitem{Bonettini-etal-2009}
S.~Bonettini, R.~Zanella, and L.~Zanni.
\newblock A scaled gradient projection method for constrained image deblurring.
\newblock {\em Inverse Probl.}, 25(1), Jan. 2009.

\bibitem{Calatroni-Chambolle-2019}
L.~Calatroni and A.~Chambolle.
\newblock Backtracking strategies for accelerated descent methods with smooth
  composite objectives.
\newblock {\em SIAM J. Optim.}, 29(3):1772--1798, 2019.

\bibitem{Chambolle-Dossal-2014}
A.~Chambolle and C.~Dossal.
\newblock On the convergence of the iterates of the "{F}ast {I}terative
  {S}hrinkage/{T}hresholding {A}lgorithm".
\newblock {\em J. Optim. Theory Appl.}, 166(3):968--982, Sept. 2015.

\bibitem{Combettes-Vu-2014}
P.~Combettes and B.~V{\~{u}}.
\newblock Variable metric forward-backward splitting with applications to
  monotone inclusions in duality.
\newblock {\em Optimization}, 63(9):1289--1318, 2014.

\bibitem{Debnath-etal-1990}
L.~Debnath and P.~Mikusi\'{n}ski.
\newblock {\em Introduction to Hilbert spaces with applications}.
\newblock Academic Press, Boston, 1990.

\bibitem{Florea2020}
M.~I. {Florea} and S.~A. {Vorobyov}.
\newblock A generalized accelerated composite gradient method: Uniting
  nesterov's fast gradient method and fista.
\newblock {\em IEEE Trans. Signal Process.}, 68:3033--3048, 2020.

\bibitem{Harmany12}
Z.~Harmany, R.~Marcia, and R.~Willett.
\newblock This is {SPIRAL-TAP}: {S}parse {P}oisson {I}ntensity {R}econstruction
  {AL}gorithms - {T}heory and practice.
\newblock {\em IEEE Trans. Image Process.}, 21(3):1084--1096, 2012.

\bibitem{Lanteri-etal-2001}
H.~Lant\'eri, M.~Roche, O.~Cuevas, and C.~Aime.
\newblock A general method to devise maximum likelihood signal restoration
  multiplicative algorithms with non-negativity constraints.
\newblock {\em Signal Process.}, 81(5), May 2001.

\bibitem{SAGE-FISTA}
S.~Rebegoldi and L.~Calatroni.
\newblock Scaled, inexact and adaptive generalized {FISTA} for strongly convex
  optimization.
\newblock {\em \url{https://arxiv.org/abs/2101.03915}}, 2021.

\bibitem{Salzo-Villa-2012}
S.~Salzo and S.~Villa.
\newblock Inexact and accelerated proximal point algorithms.
\newblock {\em J. Convex Anal.}, 19(4):1167--1192, 2012.

\bibitem{Scheinberg-2014}
K.~Scheinberg, D.~Goldfarb, and X.~Bai.
\newblock Fast first--order methods for composite convex optimization with
  backtracking.
\newblock {\em Found. Comput. Math.}, 14:389--417, 2014.

\bibitem{Schmidt2011}
M.~Schmidt, N.~L. Roux, and F.~Bach.
\newblock Convergence rates of inexact proximal-gradient methods for convex
  optimization.
\newblock {\em arXiv:1109.2415v2}, 2011.

\end{thebibliography}

\end{document}